\documentclass[11pt]{article}
\usepackage[english]{babel} %
\usepackage{amssymb,amsmath,amsthm,enumitem}
\usepackage{xcolor}
\usepackage[
    colorlinks=true,
    linkcolor=blue,   
    citecolor=blue,   
    urlcolor=blue!60 
]{hyperref}
\usepackage{xurl}
\hypersetup{pdftitle={Factorizations of finite groups},pdfauthor={Mikhail Kabenyuk}}

\theoremstyle{plain}
\newtheorem{theorem}{Theorem}[section] 
\newtheorem*{theorem*}{Theorem} 
\newtheorem{example}[theorem]{Example}
\newtheorem{proposition}[theorem]{Proposition}
\newtheorem{question}[theorem]{Question}
\newtheorem{lemma}[theorem]{Lemma}
\newtheorem{corollary}[theorem]{Corollary}

\theoremstyle{definition}
\newtheorem{remark}[theorem]{Remark} 

\def\gr{\operatorname{gr}}
\def\Aut{\operatorname{Aut}}
\def\Ker{\operatorname{Ker}}
\def\GL{\operatorname{GL}}
\def\SL{\operatorname{SL}}

\setlist[enumerate]{%
  labelwidth=!,
  labelindent=\parindent,
  leftmargin=0pt,
  labelsep=*,
  align=left,
  itemindent=\dimexpr\parindent+1.2em\relax,
  itemsep=2pt,
  topsep=2pt,
  partopsep=0pt,
  parsep=0pt,
  font=\normalfont
}

\title{Factorizations of finite groups}
\author{Mikhail Kabenyuk}
\date{}

\begin{document}\maketitle

\begin{abstract}
A finite group $G$ is called $k$-factorizable
if for every ordered factorization $|G|=a_1\cdots a_k$ into integers each greater than $1$
there exist subsets $A_1,\dots,A_k\subseteq G$ such that $|A_i|=a_i$ for each $i$ 
and $G=A_1\cdots A_k$.

The main results are as follows.

1. (Theorem 1.1)
For every integer $k\geq3$ there exists a finite group $G$
such that $G$ is not $k$-factorizable.

2. (Theorem 1.2)
Let $G$ be a finite group of order $4m$.
If a Sylow $2$-subgroup of $G$ is elementary abelian
and all involutions of $G$ are conjugate,
then $G$ has no factorization of the form $G=ABC$ with $|A|=|C|=2$ and $|B|=m$.

3. Only $8$ groups of order at most $100$ fail to be $k$-factorizable for some $k$
(Theorem 1.3).
\end{abstract}

\section{Introduction}\label{section-Introduction}
\def\Hajos{Haj$\acute{\rm{o}}$s}
\def\Redei{R$\acute{\rm{e}}$dei}
\def\Szabo{Szab$\acute{\rm{o}}$}
\def\circm{\mathbin{\vcenter{\hbox{$\scriptscriptstyle\circ$}}}}

Let $G$ be a finite group, and let $A_1,\ldots,A_k$ be subsets of $G$ such that $G=A_1\dotsm A_k$. 
If, in addition, $|G|=|A_1|\dotsm|A_k|$, 
then we call this a \textit{$k$-fold factorization of $G$ of type $(|A_1|,\ldots,|A_k|)$}, 
or, more briefly, an \textit{$(|A_1|,\ldots,|A_k|)$-factorization of $G$}.

If $G=A_1\dotsm A_k$, then the condition $|G|=|A_1|\dotsm|A_k|$ is equivalent to the requirement that 
every element $g\in G$ can be uniquely represented as $g=x_1\dotsm x_k$, 
where $x_i\in A_i$ for $i=1,\ldots,k$.
Indeed, the mapping 
$(x_1,\ldots,x_k)\mapsto x_1\dotsm x_k$ from 
$A_1\times\dotsm\times A_k$ to $G$ is surjective, and hence bijective if and only if
\[
|A_1\times\dotsm\times A_k|=|A_1|\dotsm|A_k|=|G|.
\]
To indicate that the product is a factorization, 
we sometimes write $G=A\circm B$ or $G=A_1\circm\dotsm\circm A_k$; 
however, when this is clear from the context, we often simply write $G=AB$ or $G=A_1\dotsm A_k$.

The problem of factorizing finite abelian groups as direct products of subsets arose 
in \Hajos' \cite{Hajos} solution of a classical problem of Minkowski.
Minkowski's conjecture \cite{Minkowski} states that every lattice tiling of $\mathbb{R}^n$ by unit
$n$-cubes contains at least two $n$-cubes that
meet completely on one of their $(n-1)$-dimensional faces.
In 1938, \Hajos\ reformulated Minkowski's conjecture in terms
of finite abelian groups:
In every factorization of a finite abelian group into cyclic subsets, 
at least one of the factors is a subgroup.
A subset of a group is called cyclic if it has the form $\{e,c,\ldots,c^s\}$ 
for some group element $c$ and some positive integer $s$.
In 1941, \Hajos\ proved Minkowski's conjecture in this group-theoretic form \cite{Hajos}.

In 1965, \Redei\ \cite{Redei} proved that if $G$ is a finite abelian group and 
$G=A_1\dotsm A_k$ is a factorization in which each $A_i$ contains the identity element of $G$ 
and has prime order, then at least one of the sets $A_i$ is a subgroup.
For these and other classical results on factorizations of abelian groups, 
see the book by \Szabo\ \cite{Szabo} and the book by \Szabo\ and Sands \cite{Szabo-Sands}.
Note also that in the first edition of Fuchs' book ``Abelian Groups'' 
an entire chapter is devoted to factorizations of finite abelian groups \cite[Chapter XV]{Fuchs}.
In this chapter, the results of \Hajos\ and \Redei\ are presented in detail.
In the latest edition of this book (2015), \Hajos' results are presented in Section 3 of Chapter 3.

Bernstein \cite{Bernstein} generalized \Hajos' theorem on products of 
cyclic subsets to noncommutative groups.

An important line of research over the past several decades has been the conjecture 
that every finite group $G$ admits a factorization $G=A_1\dotsm A_k$ 
in which each $A_i$ has prime order or order $4$.
Such a factorization is called a minimal logarithmic signature (MLS).
It is known that every finite solvable group has an MLS.
If $G$ is a finite group with a normal subgroup $K$ such that 
both $K$ and $G/K$ have an MLS, then $G$ has an MLS. 
The latter result implies that the existence of an MLS for an arbitrary finite group 
can be reduced to the existence of MLSs for finite simple groups \cite{Rahimipour}.

Against this background, it is natural to ask whether an arbitrary factorization of $|G|$ 
can always be realized by a factorization of the group $G$ into subsets of the corresponding sizes.

Hooshmand \cite[Problem 19.35]{Kourovka} 
(see also \cite{Hooshmand2014}, \cite[Questions I and II]{Hooshmand}, \cite{Banakh2}, and \cite{Banakh1}) 
asked the following questions:

a) \textit{%
Is it true for every finite group $G$ of order $n$
that for every factorization $n=a_1\dotsm a_k$
there exist subsets $A_1,\ldots,A_k$ of $G$
such that $|A_1|=a_1,\ldots,|A_k|=a_k$ and $G=A_1\dotsm A_k?$
}

b) \textit{%
The same question for $k=2$.
}

The case $k=2$, that is, the class of $2$-factorizable groups, 
has been studied in \cite{Bildanov} and \cite{Hooshmand}.
It was proved that every finite group is $2$-factorizable 
if and only if every finite simple group is $2$-factorizable.
All finite simple groups of order at most $10000$ are $2$-factorizable \cite{Bildanov}.

For $k=3$, the answer to Hooshmand's question (a) is negative.
Bergman \cite{Bergman} (see also \cite{Banakh1})
proved that
the alternating group $A_4$ is not $3$-factorizable since it has no $(2,3,2)$-factorization.
By computer calculation, Brunault \cite{Brunault} proved that 
the alternating group $A_5$ has no $(2,3,5,2)$-factorization.
However, this does not imply that $A_5$ is not $3$-factorizable.

Let $k\geq2$ be a fixed integer.
We say that $G$ is \textit{$k$-factorizable} if for every ordered factorization
$|G|=a_1\dotsm a_k$ into integers each greater than $1$, the group $G$
has an $(a_1,\ldots,a_k)$-factorization.
Here, ordered means that the order of the factors $a_1,\ldots,a_k$ matters.

We say that $G$ is \textit{multifold-factorizable}
if $G$ is $k$-factorizable for every integer $k\geq2$.

Our first main result shows that for every $k\geq3$ there exists
a finite group that is not $k$-factorizable.
More precisely, the following theorem holds.

\begin{theorem}\label{theorem_exist}
For every integer $k\geq3$ there exists a finite group $G$ 
and integers $a_1,\ldots,a_k>1$ such that $|G|=a_1\dotsm a_k$ and
$G$ has no $(a_1,\ldots,a_k)$-factorization.
\end{theorem}

The following theorem gives a sufficient condition for a finite group not to be $3$-factorizable.

\begin{theorem}\label{theorem_non-factorizable}
Let $G$ be a finite group of order $4m$. Assume that the following properties hold:
\begin{enumerate}
    \item[$(i)$] a Sylow $2$-subgroup of $G$ is elementary abelian;
    \item[$(ii)$] all involutions of $G$ are conjugate.
\end{enumerate}

Then $G$ has no factorization of the form $G=ABC$ with $|A|=|C|=2$ and $|B|=m$.
\end{theorem}

The next theorem gives a complete list of all
non-multifold-factorizable groups of order $\leq 100$.

\begin{theorem}\label{theorem_list_100}
There are exactly $8$ non-multifold-factorizable groups
of order at most $100$.
Here is the complete list of such groups:
\[
  \begin{array}{ll}
    1.\ A_4                        &            5.\ A_5\\
    2.\ (C_2\times C_2)\rtimes C_9 &            6.\ C_5\times A_4 \\
    3.\ C_3\times A_4              &            7.\ C_7\times A_4 \\
    4.\ (C_2\times C_2\times C_2)\rtimes C_7\hspace{90pt}&  8.\ C_7\rtimes A_4
  \end{array}
  \tag{L}\label{listL}
\]
\end{theorem}
Note that each of the eight groups in the list~\eqref{listL}
is not multifold-factorizable by Theorem~\ref{theorem_non-factorizable}.

We also recall the actions defining the three semidirect products.
For the group $G=(C_2\times C_2)\rtimes C_9$ in the list~\eqref{listL}, 
the action of $C_9=\gr(t)$ on $C_2\times C_2=\gr(a,b)$ is given by
\begin{equation}\label{formula_action_36}
    tat^{-1}=b,\ tbt^{-1}=ab.
\end{equation}
In particular, $t^3$ acts trivially on $\gr(a,b)$, so $Z(G)=\gr(t^3)\cong C_3$
and $G/Z(G)\cong A_4$.

For the group $G=(C_2\times C_2\times C_2)\rtimes C_7$ in the list~\eqref{listL},
the action of $C_7=\gr(t)$ on
$C_2\times C_2\times C_2=\gr(a,b,c)$ is given by
\begin{equation}\label{formula_action_56}
tat^{-1}=b,\ 
tbt^{-1}=c,\ 
tct^{-1}=ac.
\end{equation}
This is the unique nonabelian semidirect product of $(C_2)^3$ by $C_7$;
equivalently, it is the affine linear group $AGL(1,8)$ and can be written as
$G\cong \mathbb{F}_8\rtimes \mathbb{F}_8^*$, where $\mathbb{F}_8^*\cong C_7$
acts on the additive group $\mathbb{F}_8\cong (C_2)^3$ by multiplication.

Finally, for the group $G=C_7\rtimes A_4$ in the list~\eqref{listL},
the group $A_4=\gr(a,b,t)$ acts on $C_7=\gr(v)$ by the rules
\begin{equation}\label{formula_action_84}
av=va,\ bv=vb,\ tvt^{-1}=v^4,
\end{equation}
where $\gr(a,b)\cong C_2\times C_2$ is the normal Klein four subgroup of $A_4$.
Thus the action of $A_4$ on $C_7$ factors through
$A_4/\gr(a,b)\cong C_3$, so $G$ may also be viewed as
$(C_7\times (C_2\times C_2))\rtimes C_3$ with $C_3$ acting nontrivially 
on both $C_7$ and $C_2\times C_2$.

\medskip
We use the symbol $e$ for the identity element of a group $G$.
If $X$ is a subset of $G$, then $|X|$ denotes the cardinality of $X$, and 
$\gr(X)$ denotes the subgroup of $G$ generated by $X$.
For a singleton $X=\{x\}$, we write $\gr(x)$ for the cyclic subgroup of $G$ generated by $x$.
If $A\subseteq G$ and $x\in G$, we write $A^x=x^{-1}Ax$.
For a subset $X\subseteq G$, the symbols $C_G(X)$ and $N_G(X)$ denote 
the centralizer and the normalizer of $X$ in $G$, respectively.
If $H\leq G$, then $|G:H|$ denotes the index of $H$ in $G$.
The center of $G$ is denoted $Z(G)$.
If $X$ and $Y$ are subgroups of $G$, we denote by
$[X,Y]$ the subgroup $\gr(x^{-1}y^{-1}xy\mid x\in X,y\in Y)$.

\section{Bergman's question}\label{section-Bergman's_question}
A version of the following lemma was proved by Bergman in \cite[Lemma 2.1]{Bergman}, 
but we include a proof for completeness. 
Our argument differs only slightly from Bergman's.
See also \cite[Lemma 2.3 and Corollary 2.1]{Szabo-Sands}.

\begin{lemma}\label{lemma_divisibility}
    Let $G=A\circm B$ be a factorization of a finite group $G$. 
    If $H$ is a subgroup of $G$ containing $A$, then $|A|$ divides $|H|$.
    The same statement holds with $A$ replaced by $B$. 
\end{lemma}

\begin{proof}
 Since $A\subseteq H$, we have $AB\cap H=A(B\cap H)$, and since $G=A\circm B$, 
it follows that
\[
H=G\cap H=AB\cap H=A(B\cap H)=A\circm(B\cap H).
\]
Thus $H=A\circm(B\cap H)$ is a factorization of $H$, which implies 
$|A|$ divides $|H|$.
The corresponding statement for $B$ is proved in the same way.
\end{proof}

In this context, Bergman asked the following question \cite[Question 2.4]{Bergman}
(see also \cite[Problem 21.18]{Kourovka}).

\begin{question}
    If a finite group $G$ has a factorization
    $G=A_1\circm A_2\circm A_3$, must $|A_2|$ divide the
    order of the subgroup $H$ of $G$ generated by $A_2$?
\end{question}

The following proposition shows that the answer to this question is negative.

\begin{proposition}
    The alternating group $A_4$ admits a factorization $A_4=A\circm B\circm C$ such that
    $|B|=2$ and $|\gr(B)|=3$.
\end{proposition}

\begin{proof}
Permutations are composed from left to right, as in \cite{MHall}.
Set
\[
A=\{e,(14)(23)\},\
B=\{e,(132)\},\
C=\{e,(124),(142)\}.
\]
We show that $A_4=ABC$.
It is clear that
\[
AB=\{e,(132),(14)(23),(143)\}.
\]
Since $C$ is a subgroup of $A_4$ and $e,(132),(14)(23),(143)$
are representatives of distinct left cosets of $C$, it follows that
$A_4=ABC$.
\end{proof}

\section{Groups without \texorpdfstring{$k$}{k}-factorizations}
\label{section-Groups_without_k-factorizations}
Recall that Bergman \cite[Proposition 1.4]{Bergman}
gave a negative answer to Hooshmand's question \cite[Question 19.35]{Kourovka} for $k=3$.
What about larger $k$?
If we do not allow factorizations of $|G|$ in which some factors are equal to $1$, 
that is, if we consider only factorizations of $|G|$ into factors greater than $1$, 
then Bergman's example does not answer the question for $k>3$.
In \cite[p. 159]{Bergman}, Bergman wrote that
``With that restriction, the problem remains open for all $k > 3$;
I do not know whether there is an easy way
to modify the present example to cover those cases''.

We first prove Theorem~\ref{theorem_non-factorizable}
and then deduce Theorem~\ref{theorem_exist},
which answers this question for all $k\geq3$.

\begin{lemma}\label{lemma_normalization}
    If for some positive integers $a_1,\ldots,a_k$ the group $G$ has 
    a $k$-fold factorization $G=A_1\circm\dotsm\circm A_k$ with $|A_i|=a_i$ 
    for $i=1,\ldots,k$, then $G$ has 
    a $k$-fold factorization of the same type in which every factor contains $e$.
\end{lemma}

\begin{proof}
If $e\in A_i$ for all $i$, there is nothing to prove. 
Otherwise, let $s\geq1$ be the smallest index such that $e\notin A_s$. 
Hence $A_1,\ldots,A_{s-1}$ contain $e$. 
Choose an element $a\in A_s$. Then
\begin{align*}
  G =a^{-1}G
    & =a^{-1}A_1\circm A_2\circm\dotsm\circm A_{s-1}\circm A_s\circm\dotsm\circm A_k \\
    &=(a^{-1}A_1a)\circm(a^{-1}A_2a)\circm\dotsm\circm(a^{-1}A_{s-1}a)\circm(a^{-1}A_s)\circm\dotsm\circm A_k.
\end{align*}
For each $i<s$, the set $a^{-1}A_ia$ contains $e$, and the $s$-th factor
$a^{-1}A_s$ also contains $e$.
An easy induction completes the proof.
\end{proof}

\begin{remark}
A factorization $G=A_1\circm\dotsm\circm A_k$ is called \textit{normalized} 
if $e\in A_i$ for each $i$ (see \cite[p. 5]{Szabo-Sands}).
Thus Lemma \ref{lemma_normalization} states that 
if there is an $(a_1,\ldots,a_k)$-factorization of $G$,
then there is also a normalized $(a_1,\ldots,a_k)$-factorization of $G$.
This lemma also appears in \cite[Lemma 2.2]{Bergman}.
\end{remark}

The following lemma allows us to replace a two-element end factor
by a subgroup of order $2$ whenever its nonidentity element has order
twice an odd number.

\begin{lemma}\label{lemma_two_element_replacement}
Let $G$ be a finite group, let $a\in G$, $a\neq e$, and suppose that
\[
G=\{e,a\}\circm X
\]
is a factorization of $G$. Then the order of $a$ is even.
If $|a|=2r$ with $r$ odd, then
\[
G=\{e,a^r\}\circm X.
\]
The analogous assertions hold for a right end factor.
\end{lemma}

\begin{proof}
By uniqueness of the factorization, $G=X\sqcup aX$, where $\sqcup$
denotes disjoint union. Multiplying this partition on the left by $a$ gives
\[
G=X\sqcup aX=aX\sqcup a^2X.
\]
Both $X$ and $a^2X$ are complements of $aX$ in $G$, so $a^2X=X$.
If $a$ had odd order, then $a$ would be a power of $a^2$, and hence $aX=X$, a contradiction. 
Thus $|a|$ is even.

Now suppose that $|a|=2r$ with $r$ odd.
Since $r-1$ is even and $a^2X=X$, we have
\[
a^rX=a\cdot a^{r-1}X=aX.
\]
Thus $G=X\sqcup a^rX=\{e,a^r\}\circm X$.
The right-handed assertion follows in the same way.
\end{proof}

\begin{proof}[Proof of Theorem \normalfont\ref{theorem_non-factorizable}]
Suppose that $G=A\circm B\circm C$ is a $(2,m,2)$-factorization.
By Lemma~\ref{lemma_normalization}, we may assume that all three factors contain $e$. 
By condition~$(i)$, every element of even order in $G$ has order twice an odd number.
Apply Lemma~\ref{lemma_two_element_replacement} first to $A$ and then to $C$. 
Each replacement preserves the middle factor and the other end factor. 
We therefore obtain
\[
G=H\circm B\circm K,
\]
where $H$ and $K$ are subgroups of order $2$.

The double cosets $HbK$, $b\in B$, form a partition of $G$,
and each has $4$ elements by uniqueness of the factorization.
Since all involutions of $G$ are conjugate, there exists $x\in G$
such that $x^{-1}Hx=K$. But then
\[
HxK=x(x^{-1}Hx)K=xK,
\]
which has only $2$ elements. This contradicts the fact that every
$(H,K)$-double coset has $4$ elements.
\end{proof}

\begin{corollary}\label{corollary_omega_n}
    Let $G$ be a finite group of order $n=4m$
    satisfying the conditions of
    Theorem \textup{\ref{theorem_non-factorizable}} and
    let $\Omega(n)$ be
    the number of prime factors of $n$ (with multiplicity).
    Then
    for any integer $k$ in the interval $3\leq k\leq\Omega(n)$
    and for any integers $a_1,\ldots,a_k$ such that
    $n=a_1\dotsm a_k$ and $a_1=a_k=2$
    the group $G$ has no $(a_1,\ldots,a_k)$-factorization.
\end{corollary}

\begin{proof}
Suppose that for some $k\geq3$ the group $G$ has a $k$-fold factorization
$G=A_1\circm\dotsm\circm A_k$ with $|A_1|=|A_k|=2$.
Then $G=A_1\circm B\circm A_k$, where $B=A_2\circm\dotsm\circm A_{k-1}$. 
Since $|B|=m$, this contradicts Theorem~\ref{theorem_non-factorizable}.
\end{proof}

\begin{corollary}\label{corollary_subgroup_of_odd_index}
Suppose that $H$ is a subgroup of odd index in a finite group $G$.
If $H$ satisfies the hypotheses of Theorem~\textup{\ref{theorem_non-factorizable}},
then so does $G$. In particular, $G$ has no
$(2,|G|/4,2)$-factorization and is not $3$-factorizable.
\end{corollary}

\begin{proof}
Let $P$ be a Sylow $2$-subgroup of $H$. Since $|G:H|$ is odd, $P$ is also a Sylow $2$-subgroup of $G$. 
Hence the Sylow $2$-subgroups of $G$ are elementary abelian. 
Every involution of $G$ is conjugate in $G$ to an element of $P$. 
Since all involutions of $H$ are conjugate in $H$, all involutions of $G$ are conjugate in $G$.
The conclusion follows from Theorem~\ref{theorem_non-factorizable}.
\end{proof}

We now turn to the proof of Theorem \ref{theorem_exist}.
\begin{proof}[Proof of Theorem \normalfont\ref{theorem_exist}]
Let $k\geq3$ and put
\[
G=A_4\times C_{3^{k-3}}.
\]
The subgroup $A_4\times\{e\}$ has odd index in $G$ and satisfies
the hypotheses of Theorem~\ref{theorem_non-factorizable}.
Hence $G$ satisfies these hypotheses by
Corollary~\ref{corollary_subgroup_of_odd_index}.
Since $|G|=4\cdot3^{k-2}$ and $\Omega(|G|)=k$,
Corollary~\ref{corollary_omega_n} excludes a factorization of type
\[
(2,\underbrace{3,\ldots,3}_{k-2},2).
\]
Thus $G$ is not $k$-factorizable.
\end{proof}

\begin{example}[Affine groups]\label{example_affine}
Let $q=2^s$, where $s\geq2$, and consider the affine group
\[
G(q)=
\left(
\begin{array}{cc}
\mathbb F_q^* & \mathbb F_q\\
0 & 1
\end{array}
\right).
\]
Its translation subgroup
\[
P=
\left(
\begin{array}{cc}
1 & \mathbb F_q\\
0 & 1
\end{array}
\right)
\]
is a normal elementary abelian Sylow $2$-subgroup.
Since $G(q)/P$ has odd order $q-1$, every involution of $G(q)$ lies in $P$.
Conjugation by $\operatorname{diag}(a,1)$ sends the translation
with parameter $u$ to the translation with parameter $au$.
Thus all nonidentity elements of $P$ are conjugate in $G(q)$.

Since $|G(q)|=q(q-1)$, Theorem~\ref{theorem_non-factorizable}
excludes a factorization of type
\[
\left(2,\frac{q(q-1)}4,2\right).
\]
More generally, by Corollary~\ref{corollary_omega_n}, the group $G(q)$
is not $k$-factorizable for
\[
3\leq k\leq\Omega(q(q-1))
=s+\Omega(q-1).
\]

The group $G(q)$ is the semidirect product
$(\mathbb F_q,+)\rtimes\mathbb F_q^*$.
It is a Frobenius group with kernel $P$ and complement
$\mathbb F_q^*$.
\end{example}

\begin{remark}[Simple groups]\label{remark_simple_even}
The following finite nonabelian simple groups have elementary abelian
Sylow $2$-subgroups and a single conjugacy class of involutions
(see \cite{Walter}):
\begin{enumerate}
    \item
    $\operatorname{PSL}(2,q)$, where $q\geq4$ is a prime power
    and either $q$ is even or $q\equiv3,5\pmod8$;
    \item
    the Ree groups ${}^2G_2(3^{2s+1})$, where $s\geq1$;
    \item
    the Janko group $J_1$.
\end{enumerate}
Hence, by Theorem~\ref{theorem_non-factorizable}, each group $G$
in this list has no $(2,|G|/4,2)$-factorization.

The family $\operatorname{PSL}(2,2^s)$, $s\geq2$, together with
Corollary~\ref{corollary_omega_n}, gives another proof of
Theorem~\ref{theorem_exist}.

Note that
$\operatorname{PSL}(2,4)\cong\operatorname{PSL}(2,5)\cong A_5$,
whereas the simple groups of orders $168$ and $360$
are multifold-factorizable \cite{Kab}.
\end{remark}

\section{Supersolvable and CLT groups}
\label{section-Supersolvable_and_CLT_groups}
Let $G$ be a finite group of order $n$.
If for each divisor $m$ of $n$ the group $G$ contains a subgroup $H$ of order $m$,
then we say that $G$ has the $CLT$ property.
The term \textit{$CLT$-group} (from “converse of Lagrange's theorem”) 
was apparently first introduced by Bray \cite{Bray}.
Finite nilpotent groups provide basic examples of $CLT$-groups.
In particular, finite abelian groups and finite $p$-groups have the $CLT$ property.

Recall that a group $G$ is supersolvable if it has a finite normal series 
$\{e\}=H_0\leq H_1\leq\dotsm\leq H_r=G$ in which each factor group $H_i/H_{i-1}$ is cyclic.
Equivalently, a finite group is supersolvable if and only if it has 
a normal series each of whose factors has prime order.
Moreover, every subgroup and every factor group of a supersolvable group is supersolvable, 
and every finite nilpotent group is supersolvable (see \cite[Chap. 10]{MHall}).

The "if" part of the following lemma has a long history. 
The result was first proved by Zappa in 1940 \cite{Zappa}, 
in response to a question of Ore \cite{Ore}, who had obtained a partial result. 
McLain \cite{McLain} proved the statement in 1957,
and ten years later the same result was obtained by Deskins in \cite{Deskins}.
See also Doerk \cite{Doerk}, who gives a new proof and notes, 
with reference to McLain, that the result was already known.

In view of the above remarks, we omit the proof of the following lemma.

\begin{lemma}[On $CLT$-groups]\label{lemma_CLT_supersolvable}
  A finite group $G$ is supersolvable if and only if every subgroup of $G$, 
  including $G$ itself, has the $CLT$ property.
\end{lemma}

\begin{lemma}\label{lemma_supersolvable_multifold}
    {\rm(see also \cite[Theorem 2.1]{Hooshmand})}
    If $G$ is a finite supersolvable group,
    then $G$ is multifold-factorizable.
\end{lemma}

\begin{proof}
Let $|G|=n$ and let $n=a_1\dotsm a_k$, where $a_i>1$ for all $i$.
By Lemma \ref{lemma_CLT_supersolvable}, the group $G$ and every subgroup of $G$ have the $CLT$ property.
Using this recursively, we construct a chain of subgroups
\[
G=H_1>H_2>\dotsm>H_k>H_{k+1}=\{e\},
\]
such that $|H_i|=a_i\dotsm a_k$ for $i=1,\ldots,k$.
Then we have $|H_i:H_{i+1}|=a_i$ for every $i=1,\ldots,k$.
For each $i=1,\ldots,k$, choose a set $A_i$ of left coset representatives of $H_{i+1}$ in $H_i$. 
Then $H_i=A_iH_{i+1}$ and $|A_i|=|H_i:H_{i+1}|=a_i$.
Therefore,
\[
G=A_1H_2=A_1A_2H_3=\dotsm=A_1A_2\dotsm A_k,
\]
which is an $(a_1,\ldots,a_k)$-factorization of $G$.
\end{proof}

\begin{remark}
Note that the converse of Lemma \ref{lemma_supersolvable_multifold} does not hold in general.
For example, the symmetric groups $S_4$ and $S_5$ are not supersolvable.
However, as we shall see below, both groups are multifold-factorizable.
\end{remark}

We will need the following lemmas on supersolvable groups
in Section \ref{section-Groups_of_order_at_most_100}.

\begin{lemma}
\label{lemma_normal_sylow_abelian_quotient}
Let $P$ be a normal Sylow $p$-subgroup of a finite group $G$.
If $G/P$ is abelian and its exponent divides $p-1$, then
$G$ is supersolvable.
\end{lemma}

\begin{proof}
We proceed by induction on $|P|$. If $P=\{e\}$, then $G$
is abelian, so the assertion is clear.

Suppose that $P\neq\{e\}$, and regard
\[
W=\{x\in Z(P):x^p=e\}
\]
as a nonzero vector space over $\mathbb F_p$.
The conjugation action of $G$ on $W$ factors through $G/P$.
The corresponding linear operators commute and are annihilated
by $X^{p-1}-1$. Since this polynomial splits into distinct
linear factors over $\mathbb F_p$, these operators are
simultaneously diagonalizable. A common eigenline therefore
gives a normal subgroup $A$ of $G$ of order $p$.

The group $G/A$ satisfies the same hypotheses with the normal
Sylow $p$-subgroup $P/A$, so it is supersolvable by induction.
Lifting a normal series with cyclic factors of prime order
from $G/A$ to $G$ and adjoining $\{e\}<A$ proves that $G$
is supersolvable.
\end{proof}

\begin{corollary}\label{corollary_multifold-factorizable_2m_4m}
Let $G$ be a finite group, $p$ an odd prime, and $n$ a non-negative integer.
The group $G$ is multifold-factorizable if one of the following conditions holds:
\begin{enumerate}
    \item
    $|G|=2p^n$;
    \item
    $|G|=4p^n$, $p\equiv1\pmod4$;
    \item
    $|G|=4p^n$,
    the Sylow $p$-subgroup is normal in $G$ and
    a Sylow $2$-subgroup of $G$ is elementary abelian.
\end{enumerate}
\end{corollary}

\begin{proof}
Let $P$ be a Sylow $p$-subgroup of $G$.

In case~(1), $P$ has index $2$, so it is normal and
$G/P\cong C_2$.

In case~(2), the number of Sylow $p$-subgroups divides $4$
and is congruent to $1$ modulo $p$. Since $p\geq5$, this
number is $1$. Thus $P$ is normal, and $G/P$ is abelian
of order $4$, with exponent dividing $4\mid p-1$.

In case~(3), $P$ is normal by assumption. If $Q$ is a Sylow
$2$-subgroup of $G$, then
\[
G/P\cong Q\cong C_2\times C_2.
\]

In all three cases, $P$ is normal and $G/P$ is abelian
of exponent dividing $p-1$. Hence $G$ is supersolvable
by Lemma~\ref{lemma_normal_sylow_abelian_quotient}, and
Lemma~\ref{lemma_supersolvable_multifold} completes the proof.
\end{proof}

Every group of order $3\cdot 2^n$ with $n\geq0$ is solvable but need not be supersolvable; 
for example, $A_4$ is not supersolvable. 
However, such groups do have a useful property, which we state in the following lemma.
\begin{lemma}\label{lemma_normal_subgroup_order-2or4} 
    Let $G$ be a group of order $3\cdot 2^n$, where $n\geq2$. 
    Then $G$ contains a normal subgroup of order $2$ or $4$.
\end{lemma}

\begin{proof}
We first find a nontrivial normal $2$-subgroup $N$ such that
$G/N\cong C_3$ or $S_3$.
If a Sylow $2$-subgroup is normal, take it as $N$.
Otherwise, $G$ has exactly three Sylow $2$-subgroups.
The conjugation action on them has transitive image in $S_3$,
so its kernel $N$ is a nontrivial normal $2$-subgroup,
and $G/N\cong C_3$ or $S_3$.

Choose a minimal normal subgroup $A$ of $G$ contained in $Z(N)$.
Then $A$ is elementary abelian, and conjugation makes it
an irreducible $\mathbb F_2[G/N]$-module.

Let $T$ be the operator induced by an element of order $3$
in $G/N$. If $G/N\cong S_3$, let $C$ be the operator induced
by an involution, so that $CTC^{-1}=T^{-1}$.
If $G/N\cong C_3$, put $C=I$.
Since $(C-I)^2=0$, there is a nonzero vector $a\in A$
fixed by $C$.

The operator $C$ fixes $a$ and permutes $T(a)$ and $T^2(a)$.
Hence the nonzero subspace
\[
\langle a,T(a),T^2(a)\rangle
\]
is invariant under both $T$ and $C$, and therefore equals $A$.
Moreover, the vector
\[
b=a+T(a)+T^2(a)
\]
is fixed by both operators.
If $b\neq0$, irreducibility gives $A=\langle b\rangle$.
If $b=0$, then $A=\langle a,T(a)\rangle$.
Thus $\dim_{\mathbb F_2}A\leq2$, so $|A|=2$ or $4$.
\end{proof}

\section{The symmetric groups \texorpdfstring{$S_4$}{S4} and \texorpdfstring{$S_5$}{S5}}
\label{section-The_symmetric_groups_S_4_and_S_5}
We begin with a simple lemma.
\begin{lemma}\label{lemma_SimpleGeneral}
Let $G$ be a finite group.
\begin{enumerate}
    \item\label{inverse_order}
If $G=A_1\circm\dotsm\circm A_k$ is an $(a_1,\ldots,a_k)$-factorization of $G$, then
$G=A_k^{-1}\circm\dotsm\circm A_1^{-1}$ is an $(a_k,\ldots,a_1)$-factorization.
    \item\label{subgroup_index}
If $H$ is a subgroup of $G$, $m=|G:H|$ and $H$ has
an $(h_1,\ldots,h_k)$-factorization, then
the group $G$ possesses factorizations of the forms
\[
(h_1,\ldots,h_k,m)\quad \textit{and}\quad (m,h_1,\ldots,h_k).
\]
    \item\label{normal_subgroup}
If $H$ is a normal subgroup of $G$ and the factor group $G/H$
has an $(a_1,\ldots,a_k)$-factorization, then
the group $G$ has factorizations of each of the following forms
\[
(|H|,a_1,a_2,\ldots,a_k),\
(a_1,|H|,a_2,\ldots,a_k),\
\ldots,
(a_1,a_2,\ldots,a_k,|H|).
\]
\end{enumerate}
\end{lemma}

\begin{proof}
The first statement is obvious.
If $R$ and $S$ are systems of representatives for the right and left cosets of $H$ in $G$, 
respectively, then $|R|=|S|=m$ and $G=HR=SH$.
Hence the second statement of the lemma follows.

Let us prove the third statement.
Let $\bar{G}=G/H$ and
$\bar{G}=\bar{A}_1\circm\bar{A}_2\circm\dotsm\circm\bar{A}_k$ be
an $(a_1,\ldots,a_k)$-factorization of $\bar{G}$.
By choosing in $G$ one representative for each coset in $\bar A_i$, 
we obtain a subset $A_i\subseteq G$ such that $|A_i|=|\bar A_i|$ and
\[
G=HA_1A_2\dotsm A_k=A_1HA_2\dotsm A_k=\dots=A_1A_2\dotsm A_kH.
\]
This proves the third statement.
\end{proof}

\begin{corollary}\label{corollary_prime_index}
Let $G$ be a finite group and let $H$ be a subgroup of index $m$.
If $H$ is multifold-factorizable and $|G|=a_1\dotsm a_k$ is a factorization of $|G|$ 
such that at least one of $a_1$ and $a_k$ is divisible by $m$, 
then $G$ admits an $(a_1,\ldots,a_k)$-factorization.
\end{corollary}
\begin{proof}
Let $|G|=a_1\dotsm a_k$. By Lemma \ref{lemma_SimpleGeneral}(\ref{inverse_order}), 
we may assume that
$m=|G:H|$ divides the last factor, so that $m\mid a_k$.
Since $H$ is multifold-factorizable, $H$ has an $(a_1,\ldots,a_{k-1},a_k/m)$-factorization.
Now Lemma \ref{lemma_SimpleGeneral}(\ref{subgroup_index}) gives an
$(a_1,\ldots,a_{k-1},a_k/m,m)$-factorization of $G$. 
Replacing the last two factors by their product, 
we obtain an $(a_1,\ldots,a_k)$-factorization of $G$.
\end{proof}

\begin{corollary}\label{corollary_normal_prime_order}
Let $G$ be a finite group, and let $H$ be a normal subgroup of prime order $p$.
If the factor group $G/H$ is multifold-factorizable, then $G$ is multifold-factorizable.
\end{corollary}
\begin{proof}
Let $|G|=a_1\dotsm a_k$, and let $p=|H|$ divide $a_i$.
Since $G/H$ is multifold-factorizable, $G/H$ has an
$(a_1,\ldots,a_{i-1},a_i/p,a_{i+1},\ldots,a_k)$-factorization.
By Lemma \ref{lemma_SimpleGeneral}(\ref{normal_subgroup}), 
$G$ has a factorization with the extra factor $p$ inserted next to $a_i/p$. 
Replacing these two adjacent factors by their product, 
we obtain an $(a_1,\ldots,a_k)$-factorization of $G$.
\end{proof}

\begin{lemma}\label{lemma_groups_of_order_24}
    Every group of order $24$ is multifold-factorizable.
\end{lemma}

\begin{proof}
Let $G$ be a group of order $24$.
Since $\Omega(24)=4$,
it is enough to find $(2,2,2,3)$ and $(2,2,3,2)$ factorizations.
If $Q$ and $P$ are Sylow $2$- and $3$-subgroups of $G$ respectively,
then $G=Q\circm P$ is an $(8,3)$-factorization.
Since $Q$ has a $(2,2,2)$-factorization, 
$G$ has a $(2,2,2,3)$-factorization.

By Lemma \ref{lemma_normal_subgroup_order-2or4}
the group $G$ has a normal subgroup $A$ of order $2$ or $4$.
If $|A|=2$, then the factor group $G/A$ has a $(2,2,3)$-factorization; 
it follows from Lemma \ref{lemma_SimpleGeneral}(\ref{normal_subgroup}) that 
$G$ possesses a $(2,2,3,2)$-factorization.
If $|A|=4$, then $|G/A|=6$.
Since $G/A$ has a $(3,2)$-factorization and $A$ has a $(2,2)$-factorization,
it follows from Lemma \ref{lemma_SimpleGeneral}(\ref{normal_subgroup})
that $G$ possesses a $(2,2,3,2)$-factorization.
\end{proof}

\begin{corollary}\label{corollary_s4}
    {\rm(see also \cite[Example 2.2]{Hooshmand})}
    The symmetric group $S_4$ is multifold-factorizable.
\end{corollary}

\begin{lemma}\label{lemma_double_cosets}
    Let $G$ be a finite group and let $A$ and $B$ be subgroups 
    of orders $a$ and $b$, respectively.
    Let
    \[
    G=At_1B\cup At_2B\cup\ldots\cup At_sB
    \]
    be a decomposition of $G$ into double cosets of $A$ and $B$,
    $T=\{t_1,t_2,\ldots,t_s\}$.
    Then $G=ATB$ is an $(a,s,b)$-factorization of $G$ if and only if
    $A^x\cap B=\{e\}$ for every $x\in G$.
\end{lemma}

\begin{proof}
If $G=ATB$ is an $(a,s,b)$-factorization of $G$, then $|AxB|=ab$ for every $x\in G$.
Since $AxB=xx^{-1}AxB=xA^xB$, we have
\begin{equation}\label{coset}
|AxB|=|xA^xB|=|A^xB|=\frac{|A^x|\cdot|B|}{|A^x\cap B|}=\frac{ab}{|A^x\cap B|}.
\end{equation}
It follows that $A^x\cap B=\{e\}$ for every $x\in G$.
Conversely, if $A^x\cap B=\{e\}$ for every $x\in G$, 
then (\ref{coset}) implies that $|AxB|=ab$ for every $x\in G$.
Hence $G=ATB$ is an $(a,s,b)$-factorization of $G$.
\end{proof}

\begin{lemma}\label{lemma_s5}
    The symmetric group $S_5$ is multifold-factorizable.
\end{lemma}

\begin{proof}
Since $|S_5|=120$ and $\Omega(120)=5$,
it suffices to show that $S_5$ is $5$-factorizable.
The number of distinct ordered factorizations of $120$ into primes is equal 
to the number of permutations of the multiset
$\{2,2,2,3,5\}$, that is $5!/3!=20$.
However, by Lemma \ref{lemma_SimpleGeneral}(\ref{inverse_order})
we can omit the factorizations obtained by reversing the order of the factors.
Moreover, since $|S_5:S_4|=5$ and $S_4$ is multifold-factorizable by Corollary \ref{corollary_s4},
Corollary \ref{corollary_prime_index} shows that 
we may discard factorizations that start or end with $5$.
\begin{table}[h]
\begin{center}
\begin{tabular}{|c|l|c|c|}
  \hline
       & Factorization & $A$   & $B$\\
  \hline
  1& $(2,2,2,{\mathbf 5},3)$& $P_2$ & $P_3$ \\
  2& $(2,2,{\mathbf 5},2,3)$& $V$   & $S_3$ \\
  3& $(2,2,{\mathbf 5},3,2)$& $V$   & $S_3$ \\
  4& $(2,3,{\mathbf 2},5,2)$& $S_3$ & $D_5$ \\
  5& $(2,5,{\mathbf 2},2,3)$& $D_5$ & $S_3$ \\
  6& $(2,2,3,{\mathbf 5},2)$& $A_4=V P_3$ & $C_2$ \\
  \hline
\end{tabular}
\caption{5-factorizations of $S_5$}
\label{5-factorizations}
\end{center}
\end{table}
Therefore, we only need to consider the $6$ factorizations listed in the second column of Table \ref{5-factorizations}.

Now we use Lemma \ref{lemma_double_cosets}.
It is sufficient to choose the appropriate subgroups $A$ and $B$ of $S_5$.
These choices are specified in Table \ref{5-factorizations}, 
where we use the following notation:
$P_q$ is a Sylow $q$-subgroup of $S_5$, $q=2,3$;
$V=\{e,(12)(34),(13)(24),(14)(23)\}$;
$S_3=\gr((12),(123))$;
$A_4$ is the alternating group on four symbols 1,2,3,4;
$D_5=\gr((12345),(15)(24))$;
$C_2=\gr(12)$.
For brevity, we omit the extra pair of parentheses in expressions such as $\gr(12)$.
In all cases it is clear that $A^x\cap B=\{e\}$ for every $x\in S_5$.
The boldface number in the second column is the number of $(A,B)$-double cosets in $S_5$.
We also use suitable $2$-factorizations of the groups $V$, $S_3$, and $D_5$ to obtain the desired $5$-factorizations. 
\end{proof}

\section{More examples of factorizable groups}
\label{section-More_examples_of_factorizable_groups}
To simplify the computations, we work in the group algebra $K[G]$ over a field $K$.
For a finite subset $X\subseteq G$, define
\[
\sigma(X)=\sum_{x\in X}x.
\]
If $G=A_1\dotsm A_k$ is a factorization of $G$, then
\[
\sigma(G)=\sigma(A_1)\dotsm \sigma(A_k).
\]
Conversely, if
\[
\sigma(G)=\sigma(A_1)\dotsm \sigma(A_k),
\]
where $A_1,\dots,A_k\subseteq G$ and
\[
|A_1|\dotsm |A_k|=|G|,
\]
then $G=A_1\dotsm A_k$ is a factorization of $G$.
Indeed, for each $g\in G$, the coefficient of $g$ in the product
$\sigma(A_1)\dotsm \sigma(A_k)$ is equal in $K$ to the number of representations
$g=a_1\dotsm a_k$ with $a_i\in A_i$.
Since this coefficient is equal to $1$, each $g\in G$ has at least one such
representation. On the other hand, the total number of tuples
$(a_1,\dots,a_k)\in A_1\times\dotsm\times A_k$ is equal to
$|A_1|\dotsm |A_k|=|G|$.
Therefore each element of $G$ has exactly one representation.

If $H$ is a subgroup of $G$, then the absorption rule
\[
\sigma(H)h=h\sigma(H)=\sigma(H)
\]
holds for all $h\in H$.

\begin{lemma}\label{lemma_36-2}
    If a group $G$ has a presentation
    \begin{equation*}
        G=\gr(a,b,t\,\mid\,a^3=b^3=t^4=e,\ ab=ba,\ tat^{-1}=b,\ tbt^{-1}=a^{-1}),
    \end{equation*}
    then $G$ is multifold-factorizable.
\end{lemma}

\begin{proof}
It is clear that the group $G$ is
the semidirect product of $C_3\times C_3$ with a cyclic group $C_4$.
Since $|G|=36$ and $\Omega(36)=4$,
it suffices to show that the group $G$ is $4$-factorizable.
Since $G$ contains a subgroup $H$ of order $18$,
generated by $a$, $b$, and $t^2$,
and $H$ is supersolvable,
by Corollary \ref{corollary_prime_index}, 
it is sufficient to construct a $(3,2,2,3)$-factorization

If $A=\gr(a)$, $B=\gr(b)$, and $C=\gr(ab)$, then we have
$tAt^{-1}=t^3At^{-3}=B$ and $t^2At^{-2}=A$.
Therefore $t^kAt^{-k}\cap C=\{e\}$ for all $k=0,1,2,3$.
Let $T=\gr(t)$.
From Lemma \ref{lemma_double_cosets}
it follows that $G=A\cdot T\cdot C$ is a $(3,4,3)$-factorization.
Then $G=A\cdot\{e,t^2\}\cdot\{e,t\}\cdot C$ is
the required $(3,2,2,3)$-factorization.
\end{proof}

\begin{lemma}\label{lemma_48-1}
    If a group $G$ has a presentation
    \begin{equation*}
    G=\gr(a,b,t\,\mid\,a^4=b^4=t^3=e,\ ab=ba,tat^{-1}=a^3b^3,\ tbt^{-1}=a),
    \end{equation*}
    then $G$ is multifold-factorizable.
\end{lemma}

\begin{proof}
It is clear that the group $G$ is
the semidirect product of $C_4\times C_4$ with a cyclic group $C_3$.
Since $|G|=48$ and $\Omega(48)=5$, it is sufficient 
by Corollary \ref{corollary_prime_index} to prove 
that $G$ has factorizations of the forms $(2,2,2,3,2)$ and $(2,2,3,2,2)$.
Since the subgroup $V=\gr(a^2,b^2)$ is a normal subgroup of $G$
and the factor group $G/V$ has a $(2,2,3)$-factorization,
it follows by Lemma \ref{lemma_SimpleGeneral}(\ref{normal_subgroup})
that the group $G$ has a factorization of the form
$(2,2,3,2,2)$.

It remains to find a $(2,2,2,3,2)$-factorization.
We claim that the following equality holds in $K[G]$:
\begin{equation}\label{spec9}
    \sigma(G)=(e+a)(e+a^2)(e+t)(e+bt+a^2b^2t^2)(e+a^2b).
\end{equation}
Let $A=\gr(a)$. Then $(e+a)(e+a^2)=e+a+a^2+a^3=\sigma(A)$.
We denote the right-hand side of (\ref{spec9}) by $Q$.
Since
\[
(e+t)(e+bt+a^2b^2t^2)=(e+b^2)+(e+b)t+a(e+ab^2)t^2,
\]
and
\[
t(e+a^2b)t^{-1}=e+a^3b^2,\quad
t^2(e+a^2b)t^{-2}=e+a^3b,
\]
it follows that
\begin{equation*}
   Q=\sigma(A)\Bigl((e+b^2)(e+a^2b)+(e+b)(e+a^3b^2)t+a(e+ab^2)(e+a^3b)t^2\Bigr).
\end{equation*}
Let $B=\gr(b)$, $H=\gr(a,b)$, and $T=\gr(t)$. Then
\begin{align*}
&\sigma(A)(e+b^2)(e+a^2b) = \sigma(A)\sigma(B),\\
&\sigma(A)(e+b)(e+a^3b^2) = \sigma(A)\sigma(B),\\
&\sigma(A)a(e+ab^2)(e+a^3b) = \sigma(A)\sigma(B).
\end{align*}
Since $\sigma(A)\sigma(B)=\sigma(H)$ and $\sigma(H)\sigma(T)=\sigma(G)$, it follows that
\begin{equation*}
   Q=\sigma(A)\sigma(B)(e+t+t^2)=\sigma(H)\sigma(T)=\sigma(G). \qedhere
\end{equation*}
\end{proof}

\begin{lemma}\label{lemma_48-X}
    If a group $G$ has a presentation
    \begin{equation*}
    \begin{split}
    G=\gr(a_1,a_2,b_1,b_2,t\,\mid\,
    a_i^2=b_i^2=t^3=e,\
    a_ib_j=b_ja_i,\\
    a_1a_2=a_2a_1,\
    b_1b_2=b_2b_1,\
    ta_1t^{-1}=a_1a_2,\ ta_2t^{-1}=a_1,\\
    tb_1t^{-1}=b_1b_2,\ tb_2t^{-1}=b_1,
    i,j=1,2),
    \end{split}
    \end{equation*}
    then $G$ is multifold-factorizable.
\end{lemma}

\begin{proof}
It is easily seen that the group $G$ is
the semidirect product
$
G=V\rtimes T,
$
where $V=\gr(a_1,a_2,b_1,b_2)$ is an elementary abelian group of order 16 and
$T=\gr(t)$ is a cyclic group of order $3$.
As in Lemma \ref{lemma_48-1}, it is sufficient to prove
that the group $G$ has factorizations of the forms
$(2,2,3,2,2)$ and $(2,2,2,3,2)$.

Let $A=\gr(a_1,a_2)$ and $B=\gr(b_1,b_2)$ be subgroups of $G$.
Consider the decomposition of $G$ into double cosets of $A$ and $B$:
$G=AB\cup AtB\cup At^2B$.
Since $A\cap B=\{e\}$, $A^t\cap B=\{e\}$, and $A^{t^2}\cap B=\{e\}$, 
each of the double cosets
$AB$, $AtB$, and $At^2B$ has size $|A|\cdot|B|=16$.
Hence $G=A\cdot\{e,t,t^2\}\cdot B$ is a $(4,3,4)$-factorization.
Since each of the groups $A$ and $B$ admits a $(2,2)$-factorization, 
it follows that $G$ has a $(2,2,3,2,2)$-factorization.

Now we find a $(2,2,2,3,2)$-factorization.
Using the notation introduced at the beginning of the section, we show that
\begin{equation}\label{spec10}
\sigma(G)=(e+a_1)(e+b_1)(e+t)(e+a_2t+a_2b_1t^2)(e+a_2b_1).
\end{equation}
Let $U_1=\gr(a_1,b_1)$. Then $(e+a_1)(e+b_1)=e+a_1+b_1+a_1b_1=\sigma(U_1)$.
Let us denote the right-hand side of (\ref{spec10}) by $Q$.
Since
\[
(e+t)(e+a_2t+a_2b_1t^2)=(e+a_1b_1b_2)+(e+a_2)t+(a_1+a_2b_1)t^2
\]
and
\[
t(e+a_2b_1)t^{-1}=e+a_1b_1b_2,\quad
t(e+a_1b_1b_2)t^{-1}=e+a_1a_2b_2,
\]
we find that
\begin{align*}
   Q=\sigma(U_1)\Bigl((e+a_1b_1b_2)(e+a_2b_1)&+(e+a_2)(e+a_1b_1b_2)t\\
    &+(a_1+a_2b_1)(e+a_1a_2b_2)t^2\Bigr).
\end{align*}
Let $U_2=\gr(a_2,b_2)$. Then $\sigma(U_1)\sigma(U_2)=\sigma(V)$. Moreover,
\begin{align*}
&\sigma(U_1)(e+a_1b_1b_2)(e+a_2b_1) = \sigma(U_1)(e+b_2)(e+a_2)=\sigma(U_1)\sigma(U_2),\\
&\sigma(U_1)(e+a_2)(e+a_1b_1b_2) = \sigma(U_1)(e+a_2)(e+b_2)=\sigma(U_1)\sigma(U_2),\\
&\sigma(U_1)(a_1+a_2b_1)(e+a_1a_2b_2)
   = \sigma(U_1)(e+a_2)(e+a_2b_2)=\sigma(U_1)\sigma(U_2).
\end{align*}
Therefore
\[
Q=\sigma(U_1)\sigma(U_2)(e+t+t^2)=\sigma(V)\sigma(T)=\sigma(G). \qedhere
\]
\end{proof}

\begin{lemma}\label{lemma_72-6}
    If a group $G$ has a presentation
    \begin{equation*}
    G=\gr(a,b,t\,\mid\,a^3=b^3=t^8=e,\ ab=ba,\ tat^{-1}=ab,\ tbt^{-1}=a),
    \end{equation*}
    then $G$ is multifold-factorizable.
\end{lemma}

\begin{proof}
The group $G$ is the semidirect product
$
G=V\rtimes T,
$
where $V=\gr(a,b)$ is an elementary abelian group of order $9$ and
$T=\gr(t)$ is a cyclic group of order $8$.
Since $|G|=72$ and $\Omega(72)=5$,
it is enough to show that the group $G$ is $5$-factorizable.
Consider the subgroup $H=\gr(a,b,t^2)$ of $G$.
Since
\begin{equation*}
t^2at^{-2}=a^2b,\quad
t^2a^2bt^{-2}=a^{-1},
\end{equation*}
it follows that, writing $c=a^2b$, the group $H=\gr(a,c,t^2)$
is isomorphic to the group from Lemma \ref{lemma_36-2}
and hence it is multifold-factorizable.
Consequently, by Corollary \ref{corollary_prime_index}
it is sufficient to construct a $(3,2,2,2,3)$-factorization.

Using the notation introduced at the beginning of the section, we compute in $K[G]$.
Let $A=\gr(a)$, $B=\gr(b)$, and $C=\gr(ab)$. 
Then $\sigma(A)=e+a+a^2$, $\sigma(B)=e+b+b^2$, and $\sigma(C)=e+ab+a^2b^2$.
We prove that
\begin{equation*}
\sigma(G)=(e+a+a^2)(e+t^2)(e+t^4)(e+bt)(e+t+bt^6).
\end{equation*}
We first compute
\begin{equation}\label{spec72-6-Q}
Q=(e+t^2)(e+t^4)(e+bt)=(e+t^2+t^4+t^6)(e+bt).
\end{equation}
It is clear that
\begin{equation}\label{spec72-6-t^4-and-t^6}
  \begin{split}
    t^2bt^{-2}=ab,\quad
    t^4bt^{-4}=b^2,\quad
    t^6bt^{-6}=a^2b^2,\\
    t^3bt^{-3}=a^2b,\quad
    t^5bt^{-5}=a^2,\quad
    t^7bt^{-7}=ab^2.
  \end{split}
\end{equation}
From \eqref{spec72-6-Q} and \eqref{spec72-6-t^4-and-t^6} we obtain
\begin{align*}
Q&=e+t^2+t^4+t^6+bt+abt^3+b^2t^5+a^2b^2t^7,\\
Qt&=t+t^3+t^5+t^7+bt^2+abt^4+b^2t^6+a^2b^2,\\
Qbt^6&=bt^6+ab+b^2t^2+a^2b^2t^4+abt^7+b^2t+a^2b^2t^3+bt^5.
\end{align*}
It follows that
\begin{equation*}
R=Q(e+t+bt^6)=\sigma(C)(e+t^3+t^4+t^7)+\sigma(B)(t+t^2+t^5+t^6).
\end{equation*}
Since $\sigma(A)\sigma(B)=\sigma(V)$ and $\sigma(A)\sigma(C)=\sigma(V)$, we have 
\[
\sigma(A)R=\sigma(V)\sigma(T)=\sigma(G).\qedhere
\]
\end{proof}

\begin{lemma}\label{lemma_72-8}
If a group $G$ has a presentation
\begin{equation*}
\begin{split}
G&=\gr(a,b,x,y\mid a^3=b^3=e,\ ab=ba,\ x^4=e,\ x^2=y^2,\\
 & yxy^{-1}=x^{-1},\ xax^{-1}=a^2b,\ xbx^{-1}=ab,\ yay^{-1}=b^2,\ yby^{-1}=a),
\end{split}
\end{equation*}
then $G$ is multifold-factorizable.
\end{lemma}

\begin{proof}
The group $G$ is the semidirect product $G=V\rtimes Q_8$
of an elementary abelian group $V=\gr(a,b)$ of order $9$
with the quaternion group $Q_8=\gr(x,y)$ of order $8$.
Let $H=\gr(a,b,x)$. Then $H$ has index $2$ in $G$. 
Writing $c=a^2b$, we have $xax^{-1}=c$ and $xcx^{-1}=a^{-1}$, 
so $H=\gr(a,c,x)$ is isomorphic to the group from Lemma \ref{lemma_36-2}. 
Hence $H$ is multifold-factorizable. 
Therefore, by Corollary \ref{corollary_prime_index}, 
 it is sufficient to construct a $(3,2,2,2,3)$-factorization of $G$.
Using the notation introduced at the beginning of the section, we compute in $K[G]$.
We prove that
\begin{equation*}
\sigma(G)=(e+a+a^2)(e+bx^2)(e+b^2x)(e+y)(e+x+by).
\end{equation*}
We first compute the product
\begin{align*}
R&=(e+bx^2)(e+b^2x)(e+y)=(e+bx^2+b^2x+b^2x^3)(e+y)\\
&=
e+b^2x+bx^2+b^2x^3+
y+b^2xy+bx^2y+b^2x^3y.
\end{align*}
Each element of $Q_8$ is uniquely represented as $x^py^q$,
where $p=0,1,2,3$ and $q=0,1$.
Using $yx=x^3y$ we obtain
\begin{equation*}
Rx=
x+b^2x^2+bx^3+b^2+
x^3y+b^2y+bxy+b^2x^2y
\end{equation*}
and
\begin{equation*}
Rby=
by+axy+x^2y+a^2bx^3y+
ax^2+a^2x^3+a^2b+abx.
\end{equation*}
It follows that
\begin{align*}
S&=R(e+x+by)\\
&=
(e+b^2+a^2b)
+(b^2+e+ab)x
+(b+b^2+a)x^2\\
&+(b^2+b+a^2)x^3
+(e+b^2+b)y
+(b^2+b+a)xy\\
&+(b+b^2+e)x^2y
+(b^2+e+a^2b)x^3y.
\end{align*}
Each coefficient of $x^py^q$ in the last expression 
contains exactly one element from each coset of $A=\gr(a)$ in $V$. 
Therefore, multiplying any of these coefficients by $\sigma(A)=e+a+a^2$, 
we obtain $\sigma(V)=\sigma(A)\sigma(B)$, where $B=\gr(b)$. 
Hence 
\begin{equation*}
\sigma(A)S=\sigma(V)(e+x+x^2+x^3+y+xy+x^2y+x^3y)=\sigma(V)\sigma(Q_8)=\sigma(G).\qedhere
\end{equation*}
\end{proof}

\begin{lemma}\label{lemma_80}
    If a group $G$ has a presentation
    \begin{equation*}
    \begin{split}
    G&=\gr(a_1,a_2,a_3,a_4,t\,
    \mid\,a_i^2=t^5=e,\ a_ia_j=a_ja_i,\\
    & ta_kt^{-1}=a_ka_{k+1},
    ta_4t^{-1}=a_1,1\leq i,j\leq4,k=1,2,3),
    \end{split}
    \end{equation*}
    then $G$ is multifold-factorizable.
\end{lemma}

\begin{proof}
Let $V=\gr(a_1,a_2,a_3,a_4)$ and $T=\gr(t)$. 
Then $V$ is an elementary abelian group of order $16$, 
$T$ is a cyclic group of order $5$, and $G=V\rtimes T$.
Since $|G|=80$, $\Omega(80)=5$, 
$V$ is multifold-factorizable, and $|G:V|=5$, 
Corollary \ref{corollary_prime_index} shows 
that it is sufficient to construct factorizations of the forms $(2,2,5,2,2)$ and $(2,2,2,5,2)$.
We prove that
\begin{equation}
\sigma(G)=(e+a_1)(e+a_4)(e+a_2t+t^2+a_3t^3+a_2t^4)(e+t)(e+a_3)\label{spec80-1}
\end{equation}
and
\begin{equation}
\sigma(G)=(e+a_1)(e+a_4)(e+a_2t)(e+t+t^2+t^3+a_3t^4)(e+a_3).\label{spec80-2}
\end{equation}
To verify \eqref{spec80-1}, we compute the product of the last three factors
\begin{align*}
Q&=(e+a_2t+t^2+a_3t^3+a_2t^4)(e+t)(e+a_3)\\
&=\Bigl((e+a_2)+(e+a_2)t+(e+a_2)t^2+(e+a_3)t^3+(a_2+a_3)t^4\Bigr)(e+a_3).
\end{align*}
Using the identities
\begin{align*}
  &t(e+a_3)t^{-1}=e+a_3a_4,
  &t^2(e+a_3)t^{-2}=e+a_1a_3a_4,\\
  &t^3(e+a_3)t^{-3}=e+a_2a_3a_4,
  &t^4(e+a_3)t^{-4}=e+a_1a_2a_4
\end{align*}
we complete the computation:
\begin{align*}
Q&=
(e+a_2)(e+a_3)+
(e+a_2)(e+a_3a_4)t+
(e+a_2)(e+a_1a_3a_4)t^2\\
&+
(e+a_3)(e+a_2a_3a_4)t^3+
(a_2+a_3)(e+a_1a_2a_4)t^4.
\end{align*}
It follows that if $A=\gr(a_1,a_4)$ and $B=\gr(a_2,a_3)$, then
\begin{equation*}
(e+a_1)(e+a_4)Q=\sigma(A)\sigma(B)(e+t+t^2+t^3+t^4)=\sigma(V)\sigma(T)=\sigma(G).
\end{equation*}

Now we verify \eqref{spec80-2}.
We compute the product of the two middle factors in the same way
\begin{align*}
R&=(e+a_2t)(e+t+t^2+t^3+a_3t^4)\\
 &=((e+a_2a_3a_4)+(e+a_2)t+(e+a_2)t^2+(e+a_2)t^3+(a_3+a_2)t^4).
\end{align*}
Then we get
\begin{align*}
R(e+a_3)&=\\
&\!\!(e+a_2a_3a_4)(e+a_3)+(e+a_2)(e+a_3a_4)t+(e+a_2)(e+a_1a_3a_4)t^2\\
&\!\!+(e+a_2)(e+a_2a_3a_4)t^3+(a_3+a_2)(e+a_1a_2a_4)t^4.
\end{align*}
It follows that
\begin{equation*}
(e+a_1)(e+a_4)R(e+a_3)=\sigma(V)\sigma(T)=\sigma(G).\qedhere
\end{equation*}
\end{proof}

\section{Groups of order at most \texorpdfstring{$100$}{100}}
\label{section-Groups_of_order_at_most_100}

\begin{lemma}\label{lemma_all_3(2^n)}
Let $G$ be a group of order $2^n\cdot3$.
Then $G$ is multifold-factorizable unless $G\cong A_4$.
\end{lemma}

\begin{proof}
We consider several cases.

If $n\leq2$, then the only non-abelian group of order $6$ 
and the two non-abelian groups of order $12$ other than $A_4$ are supersolvable.
Hence every group of order $2^n\cdot3$, $n\leq2$, is multifold-factorizable unless $G\cong A_4$.

If $n=3$, then $|G|=24$, and the claim follows from Lemma \ref{lemma_groups_of_order_24}.

Assume now that $n=4$, so $|G|=48$.
By Lemma~\ref{lemma_normal_subgroup_order-2or4}, we may choose
a minimal normal subgroup $A$ of $G$ of order $2$ or $4$.
If $|A|=2$, then the factor group $G/A$ is multifold-factorizable by Lemma \ref{lemma_groups_of_order_24}.
Hence, by Corollary \ref{corollary_normal_prime_order}, the group $G$ is also multifold-factorizable.

Suppose that $|A|=4$. Then $A\cong C_2\times C_2$.
If $G/A$ is multifold-factorizable, then
Lemma~\ref{lemma_SimpleGeneral}(\ref{normal_subgroup})
yields the required factorizations of $G$, so $G$ is multifold-factorizable.
Otherwise, $G/A\cong A_4$ by the case $n\leq2$ proved above.
Let $Q$ be a Sylow $2$-subgroup of $G$ and
let $P=\gr(t)$ be a Sylow $3$-subgroup.
Since $A\leq Q$ and $G/A\cong A_4$, the subgroup $Q/A$ is the normal Klein four subgroup of $G/A$. 
Hence $Q\trianglelefteq G$. 

We claim that $Q$ is abelian.
Suppose otherwise.
By the minimality of $A$, we have $A\leq Z(Q)$ and $A=[Q,Q]$.
Choose $x,y\in Q$ such that $Q=\gr(A,x,y)$.
Since $A\leq Z(Q)$, we have $[Q,Q]=\gr([x,y])$, hence $[Q,Q]\neq A$, a contradiction.
Thus $Q$ is abelian.

Put $Q^2=\{x^2:x\in Q\}$.
Since $Q^2$ is characteristic in $Q$ and contained in $A$,
the minimality of $A$ gives
\[
Q^2=A\quad\text{or}\quad Q^2=\{e\}.
\]

Conjugation by $t$ fixes no nonidentity element of either
$A$ or $Q/A$: this follows from the minimality of $A$
and from $G/A\cong A_4$, respectively.
Consequently, $C_Q(t)=\{e\}$.
Since $Q$ is abelian, for every $x\in Q$ the product
\[
x\,(txt^{-1})\,(t^2xt^{-2})
\]
is fixed by $t$, and therefore equals $e$.

Choose $d\in Q\setminus A$ and put $c=tdt^{-1}$.
The images of $c,d$ form a basis of $Q/A$, and the preceding
identity gives
\[
tdt^{-1}=c,\quad tct^{-1}=(cd)^{-1}.
\]

Suppose first that $Q^2=A$.
The squaring map
\[
Q/A\longrightarrow A,\quad xA\longmapsto x^2,
\]
is a well-defined surjective homomorphism between groups
of order $4$, hence an isomorphism.
Thus $A=\langle c^2,d^2\rangle$, and consequently
\[
Q=\langle c,d\rangle\cong C_4\times C_4.
\]
The conjugation relations above identify $G$ with the group
in Lemma~\ref{lemma_48-1}.

Now suppose that $Q^2=\{e\}$.
Then $\langle c,d\rangle$ is a $t$-invariant complement
to $A$ in $Q$, and $tct^{-1}=cd$.
Choose $b\in A\setminus\{e\}$ and put $a=tbt^{-1}$.
Applying the same product identity to $b$, we obtain
\[
A=\langle a,b\rangle,\quad
tat^{-1}=ab,\quad tbt^{-1}=a.
\]
Together with
\[
Q=A\times\langle c,d\rangle,\quad
tct^{-1}=cd,\quad tdt^{-1}=c,
\]
these relations identify $G$ with the group
in Lemma~\ref{lemma_48-X}.

Thus in either case $G$ is multifold-factorizable.
This completes the case $|G|=48$.

Finally, let $n>4$ and assume inductively that all groups of orders
$2^{n-1}\cdot3$ and $2^{n-2}\cdot3$ are multifold-factorizable.
By Lemma \ref{lemma_normal_subgroup_order-2or4}
the group $G$ has a normal subgroup $A$ of order $2$ or $4$.
By the induction hypothesis the factor group $G/A$ is multifold-factorizable.
Hence, by Lemma \ref{lemma_SimpleGeneral}(\ref{normal_subgroup}),
the group $G$ is also multifold-factorizable.
\end{proof}

\begin{lemma}\label{lemma_all_80}
    Let $G$ be a group of order $80$.
    Then $G$ is multifold-factorizable.
\end{lemma}
\begin{proof}
If the Sylow $5$-subgroup $P$ is normal in $G$, then $G$ is a supersolvable group,
and hence, by Lemma~\ref{lemma_supersolvable_multifold}, $G$ is multifold-factorizable.

Assume that $P$ is not normal.
Then Sylow's third theorem implies that $N_G(P)=P$. 
Let $Q$ be a Sylow $2$-subgroup of $G$. 
Since $n_5=16$ and any two Sylow $5$-subgroups intersect trivially, 
the number of elements of order $5$ equals $16\cdot 4=64$. 
Hence $G$ has only $80-64-1=15$ nonidentity $2$-elements. 
Therefore all nonidentity $2$-elements of $G$ lie in $Q$, so $Q$ is the unique Sylow $2$-subgroup of $G$. 
Thus $Q\lhd G$.

Suppose that $N$ is a nontrivial normal subgroup of $G$ properly contained in $Q$. 
Then $NP$ is a subgroup of $G$, and $P$ is a Sylow $5$-subgroup of $NP$. 
Since $|N|\leq 8$, Sylow's third theorem implies that $P$ is normal in $NP$. 
Hence $N_G(P)>P$, a contradiction. 
Therefore $Q$ is a minimal normal subgroup of $G$, and hence $Q$ is elementary abelian.
Thus, if $a\in Q$, $a\neq e$, and $P=\gr(t)$,
then the subgroup $\gr(a,tat^{-1},t^2at^{-2},t^3at^{-3})$ is 
a nontrivial $P$-invariant subgroup of $Q$, 
hence normal in $G$. 
Since $Q$ is a minimal normal subgroup of $G$, 
it follows that $\gr(a,tat^{-1},t^2at^{-2},t^3at^{-3})=Q$.
Set
\[
a_1=a,\quad
a_2=a\cdot tat^{-1},\quad
a_3=a\cdot t^2at^{-2},\quad
a_4=a\cdot tat^{-1}\cdot t^2at^{-2}\cdot t^3at^{-3}.
\]
It is easy to check that
\[
ta_1t^{-1}=a_1a_2,\quad
ta_2t^{-1}=a_2a_3,\quad
ta_3t^{-1}=a_3a_4,\quad
ta_4t^{-1}=a_1.
\]
The last equality follows from
\[
a\cdot tat^{-1}\cdot t^2at^{-2}\cdot t^3at^{-3}\cdot t^4at^{-4}=
(at)^5t^{-5}=e.
\]
We see that $G$ is isomorphic to the group
described in Lemma \ref{lemma_80}, and hence it is multifold-factorizable.
\end{proof}

\begin{lemma}\label{lemma_all_56}
    Let $G$ be a group of order $56$.
    Then $G$ is multifold-factorizable unless 
    $G\cong(C_2\times C_2\times C_2)\rtimes C_7$ from the list~\eqref{listL}.
\end{lemma}
\begin{proof}
By the same argument as in Lemma~\ref{lemma_all_80}, 
we may reduce to the case where the Sylow $2$-subgroup $Q$ is a minimal normal subgroup of $G$
and hence is elementary abelian.
However, in the present case we have $|Q|=8$, 
so $G$ is isomorphic to the group $(C_2\times C_2\times C_2)\rtimes C_7$ from the list~\eqref{listL}.
The action of $C_7$ on $C_2\times C_2\times C_2$ is determined by \eqref{formula_action_56}.
\end{proof}

\begin{lemma}\label{lemma_all_36}
    Let $G$ be a group of order $36$.
    Then $G$ is multifold-factorizable 
    unless $G\cong C_3\times A_4$ or $G\cong (C_2\times C_2)\rtimes C_9$, 
    both from the list~\eqref{listL}.
\end{lemma}
\begin{proof}
If $G$ is supersolvable, then $G$ is multifold-factorizable by Lemma~\ref{lemma_supersolvable_multifold}.
Therefore, we may assume that $G$ is not supersolvable.

Let $P$ be a Sylow $3$-subgroup of $G$, and
assume that $P\trianglelefteq G$.
Then $P$ cannot be cyclic, since otherwise $G$ would be supersolvable.
Hence $P\cong C_3\times C_3$.

If $P$ contained a subgroup of order $3$ that is normal in $G$, 
then $G$ would be supersolvable, which is a contradiction.
Hence $P$ is a minimal normal subgroup of $G$.

Let $Q$ be a Sylow $2$-subgroup of $G$.
If $Q$ is elementary abelian, then
Lemma~\ref{lemma_normal_sylow_abelian_quotient}
shows that $G$ is supersolvable, a contradiction.
Thus $Q$ is cyclic of order $4$.
Let $Q=\gr(t)$.

The subgroup $P$ may be viewed as a $2$-dimensional vector space over $\mathbb{F}_3$.
Conjugation by $t$ induces a linear automorphism $\widehat t$ of $P$.
Since $t^4=e$, we have $\widehat t^{\,4}=\mathrm{id}$.

Since $P$ is  a minimal normal subgroup of $G$, 
the automorphism $\widehat t$ has no nontrivial invariant subspaces.
Hence the characteristic polynomial of $\widehat t$ is irreducible over $\mathbb{F}_3$.
But the characteristic polynomial of $\widehat t$ divides $\lambda^4-1$, and over $\mathbb{F}_3$ we have
\[
\lambda^4-1=(\lambda-1)(\lambda+1)(\lambda^2+1).
\]
Therefore the characteristic polynomial of $\widehat t$ is $\lambda^2+1$.

Choose $a\in P$, $a\neq e$, and set $b=tat^{-1}$.
Then $a$ and $b$ are linearly independent, since otherwise $\langle a\rangle$ would be $\widehat t$-invariant.
Thus $a$ and $b$ form a basis of $P$.
Since $\widehat t^2=-\mathrm{id}$, we obtain $\widehat t(b)=a^{-1}$.
Therefore
\[
tat^{-1}=b,
\quad
tbt^{-1}=a^{-1}.
\]
We see that $G$ is isomorphic to the group described in Lemma \ref{lemma_36-2}.
Thus, if a Sylow $3$-subgroup of $G$ is normal, then $G$ is multifold-factorizable.

\medskip
Assume now that $P$ is not normal in $G$.
By Sylow's theorem, the number $n_3$ of Sylow $3$-subgroups is equal to $4$.
Hence the conjugation action of $G$ on the set of Sylow $3$-subgroups yields a transitive homomorphism
\[
f\colon G\to S_4.
\]
Since $f(G)$ is a transitive subgroup of $S_4$ and $|G|=36$, 
it follows that $|f(G)|=4$ or $12$.
If $|f(G)|=4$, then $|\Ker(f)|=9$, so $\Ker(f)$ is a normal Sylow $3$-subgroup of $G$,
contrary to our assumption.
Therefore $|f(G)|=12$, so $\Ker(f)$ is a normal subgroup $S$ of order $3$.
Since $G/S\cong f(G)$ has order $12$ and $G$ is not supersolvable, it follows that
\[
G/S\cong A_4,
\]
because all groups of order $12$ other than $A_4$ are supersolvable.

Let $s\in S$, $s\neq e$.
Since $\Aut(S)$ has order $2$, we have $|G:C_G(s)|\leq2$.
If $|G:C_G(s)|=2$, then
\[
|C_G(s)/S|=6,
\]
which is impossible, since $A_4$ has no subgroups of order $6$.
Thus $S\leq Z(G)$.

Let $Q$ be a Sylow $2$-subgroup of $G$.
Then $Q\cong C_2\times C_2$ and $SQ\trianglelefteq G$.
Choose $t\in G\setminus SQ$ and set $T=\gr(t)$.
Assume first that $|t|=3$.
Then
\[
G=(SQ)T=S(QT).
\]
Since $t\notin SQ$, we have $T\cap S=\{e\}$, and therefore $QT\cap S=\{e\}$.
Hence
\[
QT\cong G/S\cong A_4,
\]
so $G\cong C_3\times A_4$, which appears in the list~\eqref{listL}.

Assume now that $|t|=9$.
Then $S\leq T$, and therefore
\[
G=QT.
\]
Thus $G$ is a semidirect product of $Q$ by $T\cong C_9$,
with the action given by \eqref{formula_action_36}.
Hence $G\cong (C_2\times C_2)\rtimes C_9$, which appears in the list~\eqref{listL}.
\end{proof}

\begin{remark}\label{remark_except_232_and_2332}
Although the groups $(C_2\times C_2)\rtimes C_9$ and $C_3\times A_4$ are not multifold-factorizable, 
they admit every $4$-factorization except for a factorization of the form $(2,3,3,2)$.
This follows from Lemma \ref{lemma_SimpleGeneral}(\ref{normal_subgroup}),
since both these groups have a normal subgroup of order $3$,
and all groups of order $12$ have factorizations of the forms $(2,2,3)$ and $(3,2,2)$.
We will need this remark later on.
\end{remark}

\begin{lemma}\label{lemma_all_72}
    Every group of order $72$ is multifold-factorizable.  
\end{lemma}

\begin{proof}
Since $\Omega(72)=5$, arguing as in Lemma~\ref{lemma_s5}, 
it is enough to construct the following factorizations:
\begin{equation}
\begin{split}
(2,2,2,3,3), (2,2,3,2,3),(2,2,3,3,2),\\
(2,3,2,3,2),(2,3,2,2,3),(3,2,2,2,3).
\end{split}
\label{factorization-22233}
\end{equation}

Let $G$ be a group of order $72$. We consider several cases.

\textsc{Case 1.}
If $G$ has a normal subgroup of order $2$ or $3$, then $G$ is multifold-factorizable.

Indeed, if $G$ contains a normal subgroup $H$ of order $3$,
then the factor group $G/H$ has order $24$ and hence is multifold-factorizable
by Lemma~\ref{lemma_groups_of_order_24}.
Then $G$ is multifold-factorizable by Corollary~\ref{corollary_normal_prime_order}.

If $G$ contains a normal subgroup $H$ of order $2$,
then the factor group $G/H$ is of order $36$ and
by Remark~\ref{remark_except_232_and_2332}
it has factorizations of the forms:
$(2,2,3,3)$, $(2,3,2,3)$, and $(3,2,2,3)$.
It follows from Lemma~\ref{lemma_SimpleGeneral}(\ref{normal_subgroup})
that $G$ has a factorization of each form listed in \eqref{factorization-22233}.

\textsc{Case 2.}
If $G$ contains a subgroup of order $36$, 
then $G$ admits every $5$-factorization listed in \eqref{factorization-22233}, 
except possibly one of the form $(3,2,2,2,3)$.
This follows from Lemma~\ref{lemma_SimpleGeneral}(\ref{inverse_order}, \ref{subgroup_index})
and Remark~\ref{remark_except_232_and_2332}.

\textsc{Case 3.}
Assume that $G$ has four Sylow $3$-subgroups.
Then the conjugation action of $G$ on these subgroups induces a transitive homomorphism
\[
f:\,G\to S_4.
\]
Hence $|f(G)|\in\{4,8,12,24\}$.
If $|f(G)|=4$ or $8$, then $|\Ker(f)|=18$ or $9$, 
so $\Ker(f)$ contains a characteristic subgroup of order $9$, 
which is therefore normal in $G$, a contradiction.
If $|f(G)|=12$ or $24$, then $|\Ker(f)|=6$ or $3$, 
and therefore $\Ker(f)$ contains a characteristic subgroup of order $3$.
Thus $G$ contains a normal subgroup of order $3$, and Case~1 applies.

\textsc{Case 4.}
Assume that $G$ has a unique Sylow $3$-subgroup $P$.
If $P$ is cyclic, then $P^3$ is a normal subgroup of order $3$, so Case~1 applies.
Hence we may assume that $P\cong C_3\times C_3$.
Let $Q$ be a Sylow $2$-subgroup of $G$.
The conjugation action of $Q$ on $P$ induces a homomorphism
\[
f:\,Q\to \Aut(P).
\]
If $\Ker(f)\neq\{e\}$, then $\Ker(f)\cap Z(Q)\neq\{e\}$.
Let $L\leq \Ker(f)\cap Z(Q)$ be a subgroup of order $2$.
Since $L$ centralizes $P$ and is normal in $Q$, we have $L\trianglelefteq G$.
Therefore $G$ is multifold-factorizable by Case~1.
Thus we may assume that $\Ker(f)=\{e\}$ and identify $Q$ with a subgroup of $\Aut(P)$.

Since $P\cong C_3\times C_3$, we may regard $P$ as a $2$-dimensional vector space over $\mathbb{F}_3$.
Thus $G\cong P\rtimes L$ for some subgroup
$L\leq\Aut(P)\cong\GL(2,3)$ of order $8$.
Choose a Sylow $2$-subgroup $K$ of $\GL(2,3)$ containing $L$.

A Sylow $2$-subgroup of $\GL(2,3)$ is isomorphic to the semidihedral group
$SD_{16}$, so we may write
\[
K=\gr(r,s\mid r^8=s^2=e,\ srs=r^3).
\]
Choose generators $a,b$ of $P$ so that
\[
 r(a)=ab,\quad r(b)=a,\quad s(a)=a^{-1},\quad s(b)=ab.
\]

Since $[K,K]=\gr(r^2)$, we have $K/[K,K]\cong C_2\times C_2$.
Hence $K$ has exactly three subgroups of index $2$, namely
\[
K_1=\gr(r)\cong C_8,\quad
K_2=\gr(r^2,s)\cong D_4,\quad
K_3=\gr(r^2,rs)\cong Q_8.
\]
From now on, we identify $P$ and each $K_i$ with their canonical images in 
$G_i=P\rtimes K_i$, $i=1,2,3$.

In the group $G_1$, let $t=r$.
Then
\[
a^3=b^3=t^8=e,\quad ab=ba,\quad tat^{-1}=ab,\quad tbt^{-1}=a.
\]
Thus $G_1$ is multifold-factorizable by Lemma~\ref{lemma_72-6}.
Now let $x=r^2$ and $y=r^{-1}s$ in $G_3$.
Then all the relations from Lemma~\ref{lemma_72-8} hold in $G_3$,
and hence $G_3$ is multifold-factorizable.
It remains to consider $G_2$.
By Case~2, it suffices to construct a $(3,2,2,2,3)$-factorization.
Let
\[
t_1=r^2s,\quad t_2=r^6s,\quad u=s,
\]
and
\[
T=\gr(t_1,t_2),\quad S=\gr(u),\quad V=\gr(a,b),\quad A=\gr(a),\quad B=\gr(b).
\]
Since $V$ is a normal subgroup of $G_2$, we have $VT=TV$.
Also, $K_2=TS$, and from the action of $s$, $r^2s$, and $r^6s$ on $P$ we obtain
$SA=AS$ and $TB=BT$.
Therefore
\[
G_2=VK_2=VTS=TVS=TBAS=BTSA.
\]
Since $|T|=4$, $G_2$ has a $(3,2,2,2,3)$-factorization.
Thus every group of order $72$ is multifold-factorizable.
\end{proof}

\begin{lemma}\label{lemma_all_75}
    Every group of order $75$ is multifold-factorizable.  
\end{lemma}
\begin{proof}
Let $G$ be a group of order $75$.
By Lemma~\ref{lemma_supersolvable_multifold} we may assume that $G$ is not supersolvable.
Let $P$ be a Sylow $5$-subgroup of $G$, and let $Q$ be a Sylow $3$-subgroup of $G$.
By Sylow's theorem, $P$ is normal in $G$. 
Since $G$ is not supersolvable, $P$ cannot be cyclic. Hence $P\cong C_5\times C_5$.
Since $P\cong C_5\times C_5$ and $G=PQ$, $G$ has a $(5,5,3)$-factorization.
The group $P$ contains six cyclic subgroups of order $5$, 
so there exist two such subgroups $A$ and $B$ that are not conjugate in $G$.
It follows from Lemma~\ref{lemma_double_cosets} that $G$ has a $(5,3,5)$-factorization.
Thus every group of order $75$ is multifold-factorizable.
\end{proof}

\begin{lemma}\label{lemma_all_60}
    Let $G$ be a group of order $60$.
    Then $G$ is multifold-factorizable 
    unless $G\cong C_5\times A_4$ or $G\cong A_5$, 
    both from the list~\eqref{listL}.
\end{lemma}
\begin{proof}
Let $G$ be a group of order $60$.
By Lemma~\ref{lemma_supersolvable_multifold} we may assume that $G$ is not supersolvable.

Let $P$ be a Sylow $5$-subgroup of $G$.
It is well known that if $P$ is not normal in $G$, then $G\cong A_5$.

Assume that $P$ is normal in $G$.
Then a Sylow $3$-subgroup is not normal in $G$
(otherwise $G$ would be supersolvable), and hence $G/P\cong A_4$.
The conjugation action yields a homomorphism $f\colon G\to C_4$,
since $\Aut(P)\cong C_4$.
Hence $|\Ker(f)|=15$, $30$, or $60$.
The first two possibilities are impossible,
because then $\Ker(f)/P$ would be a normal subgroup of $G/P\cong A_4$ of order $3$ or $6$.
Consequently, $P\leq Z(G)$. 
Let $Q$ and $R$ be Sylow $2$- and $3$-subgroups of $G$, respectively. 
Since $PQ/P$ is a normal Sylow $2$-subgroup of $G/P\cong A_4$, the subgroup $PQ$ is normal in $G$. 
Since $P\leq Z(G)$, we have $Q\trianglelefteq PQ$, and therefore $Q\trianglelefteq G$. 
Thus $QR$ is a subgroup of $G$ of order $12$, and $G=P\times QR$. 
Since $G/P\cong A_4$, we have $QR\cong A_4$. 
Therefore $G\cong C_5\times A_4$.
\end{proof}

\begin{lemma}\label{lemma_all_84}
    Let $G$ be a group of order $84$.
    Then $G$ is multifold-factorizable 
    unless $G\cong C_7\times A_4$ or $G\cong C_7\rtimes A_4$, 
    both from the list~\eqref{listL}.
\end{lemma}
\begin{proof}
Let $G$ be a group of order $84$.
By Lemma~\ref{lemma_supersolvable_multifold} we may assume that $G$ is not supersolvable.

Let $P$ be a Sylow $7$-subgroup of $G$.
By Sylow's theorem, $P$ is normal in $G$ and $G/P\cong A_4$.
Then $G$ acts on $P$ by conjugation.
Since $\Aut(P)\cong C_6$, this action yields a homomorphism $f\colon G\to C_6$.
Hence $|\Ker(f)|=14$, $28$, $42$, or $84$, and therefore $|\Ker(f)/P|=2$, $4$, $6$, or $12$.
Since $A_4$ has no normal subgroups of orders $2$ and $6$,
it follows that either $|\Ker(f)/P|=12$, or $|\Ker(f)/P|=4$.

Let $Q$ be a Sylow $2$-subgroup of $G$ and $R$ a Sylow $3$-subgroup of $G$.
If $|\Ker(f)/P|=4$, then $\Ker(f)=PQ=P\times Q$.
Since $Q\trianglelefteq G$, the product $QR$ is a subgroup of $G$, and
\[
G=(P\times Q)\rtimes R.
\]
Moreover, $R$ acts nontrivially on both $P$ and $Q$, so
\[
G\cong C_7\rtimes A_4.
\]
If $|\Ker(f)/P|=12$, then $\Ker(f)=G$, so $P\leq Z(G)$.
Also, $PQ/P$ is the normal Klein four subgroup of $G/P\cong A_4$.
Then $PQ=P\times Q$ is a normal subgroup of $G$, and hence $Q\trianglelefteq G$.
Therefore the product $QR$ is a subgroup of $G$, and
\[
G=(P\times Q)\rtimes R.
\]
Since the action of $R$ on $P$ is trivial, we get
\[
G\cong C_7\times A_4.
\]
Thus, in both cases, we obtain the groups from the list~\eqref{listL}.
\end{proof}

\begin{proof}[Proof of Theorem \normalfont\ref{theorem_list_100}]
We now complete the proof of the theorem.
Since, as noted in the introduction, 
each group in the list~\eqref{listL} is not multifold-factorizable 
by Theorem~\ref{theorem_non-factorizable}, 
it remains to show that every group outside this list is multifold-factorizable.

Let $G$ be a finite group of order $|G|\leq 100$.
By Lemma~\ref{lemma_supersolvable_multifold} we may assume that $G$ is not supersolvable.

Assume first that the order $|G|$ has at least three distinct prime divisors.
Since $|G|\leq 100$, we have
\[
|G|\in\{30,42,60,66,70,78,84,90\}.
\]
If $|G|\in\{30,42,66,70,78\}$, then $|G|$ is squarefree, so $G$ is supersolvable, a contradiction.
If $|G|=60$, then by Lemma~\ref{lemma_all_60}
$G$ is multifold-factorizable or $G\cong A_5$ or $G\cong C_5\times A_4$,
both of which appear in the list~\eqref{listL}.
If $|G|=84$, then by Lemma~\ref{lemma_all_84}
$G$ is multifold-factorizable or $G\cong C_7\times A_4$ or $G\cong C_7\rtimes A_4$,
both from the list~\eqref{listL}; in the latter group the action of $A_4$ on $C_7$
is given by \eqref{formula_action_84}.
If $|G|=90$, then a Sylow $3$-subgroup of $G$ is normal; this is easy to verify.
Hence $G$ has a normal subgroup of order $45$.
Since every group of order $45$ is abelian, it follows that $G$ is supersolvable, a contradiction.

Assume next that
\[
|G|=p^r q^s,
\]
where $p<q$ are odd primes and $r,s\geq 1$.
If $r=s=1$, then $|G|$ is squarefree, so $G$ is supersolvable, a contradiction.
If $r=2$, then necessarily $s=1$; since $|G|\leq 100$, we have
\[
|G|=9q,\quad q\in\{5,7,11\}.
\]
In each of these cases, a Sylow $q$-subgroup is normal, and hence $G$ is supersolvable, a contradiction.
If $s=2$, then necessarily $r=1$; since $|G|\leq 100$, we have $|G|=75$.
By Lemma~\ref{lemma_all_75}, every group of order $75$ is multifold-factorizable.

Assume next that
\[
|G|=2^r p^s,
\]
where $p$ is an odd prime and $s>1$.
Since $|G|\leq 100$, we have
\[
|G|\in\{18,36,50,54,72,98,100\}.
\]
If $|G|\in\{18,50,54,98,100\}$, then $G$ is supersolvable.
If $|G|=36$, then by Lemma~\ref{lemma_all_36}
$G$ is multifold-factorizable
or $G\cong C_3\times A_4$
or $G\cong (C_2\times C_2)\rtimes C_9$ from the list~\eqref{listL},
where in the latter group the action of $C_9$ on $C_2\times C_2$ is given by \eqref{formula_action_36}.
If $|G|=72$, then by Lemma~\ref{lemma_all_72}
$G$ is multifold-factorizable.

Finally, assume that
\[
|G|=2^r p,
\]
where $p$ is an odd prime.
We claim that in this case $G$ is multifold-factorizable unless
$G\cong A_4$ or
$G\cong(C_2\times C_2\times C_2)\rtimes C_7$.
Both exceptions occur in the list~\eqref{listL}.

Let $P$ be a Sylow $p$-subgroup of $G$.
If $P$ is normal in $G$, then $G$ is supersolvable, since $P$ is cyclic and $G/P$ is nilpotent.

Assume that $P$ is not normal in $G$.
Let $n_p$ be the number of Sylow $p$-subgroups of $G$.
Then
\[
n_p>1,\quad n_p\equiv 1 \pmod p,\quad n_p\mid 2^r.
\]
Hence
\[
|G|=2^r p\geq n_p p>p^2.
\]
Since $|G|\leq 100$, it follows that $p\in\{3,5,7\}$.

If $|G|=2^r\cdot 3$, then by Lemma~\ref{lemma_all_3(2^n)}
$G$ is multifold-factorizable or $G\cong A_4$ from the list~\eqref{listL}.

If $|G|=2^r\cdot 5\leq 100$ and a Sylow $5$-subgroup of $G$ is not normal, then $|G|=80$.
By Lemma~\ref{lemma_all_80}, $G$ is multifold-factorizable.

If $|G|=2^r\cdot 7\leq 100$ and a Sylow $7$-subgroup of $G$ is not normal, then $|G|=56$.
By Lemma~\ref{lemma_all_56},
$G$ is multifold-factorizable
or $G\cong (C_2\times C_2\times C_2)\rtimes C_7$ from the list~\eqref{listL},
where the action of $C_7$ on $C_2\times C_2\times C_2$ is given by \eqref{formula_action_56}.

We have considered all cases, and the proof of Theorem~\ref{theorem_list_100} is complete.
\end{proof}

\section{A question}
\label{section-Questions}
Theorem~\ref{theorem_non-factorizable} gives an obstruction involving
two-element end factors. An analogous obstruction for odd primes
would require a different argument.

\medskip\textsc{Question.}
Let $H=\SL(2,p)$, let $S$ be a Sylow $2$-subgroup of $H$,
and let $V=\mathbb F_p^2$ be its natural module.
Let $G$ be one of the following groups:
\begin{enumerate}
    \item $V\rtimes S$, if $p=7$;
    \item $V\rtimes N_H(S)$, if $p=5$ or $p=11$.
\end{enumerate}
Is it true that $G$ admits no factorization of type
\[
(p,|G|/p^2,p)?
\]
In all three cases the acting group is transitive on the
one-dimensional subspaces of $V$.

For comparison, when $p=3$, the group $V\rtimes S$ has order $72$.
By Lemma~\ref{lemma_all_72}, it is multifold-factorizable and, in particular,
admits a factorization of type $(3,8,3)$.
Thus the case $p=3$ is excluded from the question.

\end{document}